\documentclass[10pt,reqno]{amsart}
\usepackage{amsmath,latexsym,epsfig,amsfonts,amssymb,graphicx}

\newtheorem{theorem}{Theorem}[section]
\newtheorem{lemma}[theorem]{Lemma}
\newtheorem{proposition}[theorem]{Proposition}

\theoremstyle{definition}

\newtheorem{opening}[theorem]{Question}
\theoremstyle{plain}

\newcommand{\R}{\ensuremath{\mathbb{R}}}

\newcommand{\PLF}{\textup{PLF}}
\newcommand{\PL}{\textup{PL}}

\newcommand{\PLc}{\textup{PLc}}
\newcommand{\PLE}{\textup{PLE}}
\newcommand{\Ho}{\ensuremath{\textup{H}(\mathbb{R})} }
\newcommand{\Hop}{\ensuremath{\textup{H}^+(\mathbb{R})}} 
\newcommand{\Hon}{\ensuremath{\textup{H}^-(\mathbb{R})} }

\numberwithin{equation}{section}

\renewcommand{\leq}{\leqslant}
\renewcommand{\geq}{\geqslant}

\begin{document}

\title[Reversibility of piecewise linear homeomorphisms]{Reversible maps and composites of involutions in groups of piecewise linear homeomorphisms of the real line}\thanks{The authors thank Anthony O'Farrell and the referees for their detailed and helpful remarks. The second author's work was supported by Science Foundation Ireland grant 05/RFP/MAT0003.}

\author{Nick Gill}
\email{nickgill@cantab.net}
\address{Institute of Mathematical Sciences, CIT Campus, Taramani, Chennai 600113, India}

\author{Ian Short}
\email{ian.short@nuim.ie}
\address{Mathematics Department, Logic House, N.U.I. Maynooth, Maynooth, County Kildare, Ireland}



\date{\today}


\begin{abstract}
An element of a group is \emph{reversible} if it is conjugate to its own inverse, and it is \emph{strongly reversible} if it is conjugate to its inverse by an involution. A group element is strongly reversible if and only if it can be expressed as a composite of two involutions. In this paper the reversible maps, the strongly reversible maps, and those maps that can be expressed as a composite of involutions are determined in certain groups of piecewise linear homeomorphisms of the real line.
\end{abstract}

\maketitle

\section{Introduction}\label{S: introduction}


A member $g$ of a group $G$ is \emph{reversible} if it is conjugate to its own inverse within $G$, and $g$ is \emph{strongly reversible} if it is reversible and the conjugating map can be chosen to be an involution. Note that $g$ is strongly reversible if and only if it can be expressed as a composite of two involutions. The object of this paper is to study reversible maps, strongly reversible maps, and composites of involutions and reversible maps in certain groups of piecewise linear homeomorphisms of the real line.

Interest in reversibility originates from the theory of dynamical systems. For example, Birkhoff (\cite{Bi15}) used time-reversal symmetry in his study of the three body problem, and Arnol'd and Sevryuk (\cite{ArSe86}) identified reversible dynamical systems that are not strongly reversible. (See \cite{LaRo98} for many more references.) 


Recently, Jarczyk and Young (\cite{Ja02b,Yo94}) classified the strongly reversible maps in the group of homeomorphisms of the real line. Related results about composites of involutions in the group of homeomorphisms of the real line can be found in \cite{FiSc55,Ja02a,McSt85,OF04}. In this paper, we prove similar results for groups of piecewise linear homeomorphisms of the real line.

\subsection{Definitions}\label{SS: formal}\quad

Let $G$ be a group. An element $g$ in $G$ is \emph{reversible} if and only if there is an element $h$ in $G$ such that $hgh^{-1}=g^{-1}$.  We say that $h$ \emph{reverses} $g$. Note that the equation $hgh^{-1}=g^{-1}$ can be manipulated to yield $h^{-1}gh=g^{-1}$. Furthermore, if $n$ is an integer then 
\begin{equation}\label{E: generalReverse}
h^ngh^{-n}=g^{(-1)^n}.
\end{equation} 
An \emph{involution} in $G$ is an element whose square is the identity. An element $g$ in $G$ is \emph{strongly reversible} if and only if there is an involution $\sigma$ such that $\sigma g\sigma^{-1}=g^{-1}$.  Define
\[
I_n(G) = \{\tau_1\tau_2\dotsb\tau_n\,|\,\text{$\tau_1,\tau_2,\dots,\tau_n$ are involutions in $G$}\,\},
\]
for $n=1,2,\dotsc$. Then
\[
I_1(G)\subseteq I_2(G) \subseteq I_3(G)\subseteq \dotsb.
\]
The set $I_1(G)$ consists of the involutions in $G$, and $I_2(G)$ consists of the strongly reversible members of $G$. The union of the sets $I_n(G)$ is a normal subgroup denoted $I_\infty(G)$; it is the subgroup of $G$ generated by involutions. Likewise we define
\[
R_n(G) = \{g_1g_2\dotsb g_n\,|\,\text{$g_1,g_2,\dots,g_n$ are reversible elements of $G$}\,\},
\]
for $n=1,2,\dotsc$. Again we have a chain of subsets of $G$,
\[
 R_1(G)\subseteq R_2(G) \subseteq R_3(G)\subseteq \dotsb.
\]
Membership of $I_n(G)$ and $R_n(G)$ is preserved under conjugation. 

We will examine the sets $R_n(G)$ and $I_n(G)$ for the groups $G$ described in the next section. 

\subsection{The groups of interest}\quad\label{SS: groups of interest}

We consider reversibility in three groups.

\begin{enumerate}
\item The full group of homeomorphisms of the real line, denoted $\textup{H}(\R)$.

\item The group $\PLF(\R)$ of piecewise linear homeomorphisms of $\R$ which are locally affine at all but a finite number of points of $\R$.

\item The group $\PL(\R)$ of piecewise linear homeomorphisms of $\R$ which are locally affine at all but a discrete set of points of $\R$.
\end{enumerate}


Each of $\Ho$, $\PLF(\R)$, and $\PL(\R)$ has a subgroup of index $2$ consisting of  orientation preserving (strictly increasing) maps. We denote these subgroups by $\Hop$, $\PLF^+(\R)$, and $\PL^+(\R)$; the corresponding cosets of orientation reversing maps are denoted by $\Hon$, $\PLF^-(\R)$, and $\PL^-(\R)$.

A word of warning: the property of being reversible (and similarly the property of being a member of one of the classes $I_n(G)$ and $R_n(G)$) depends on the group. Thus, for example, an orientation preserving homeomorphism may be reversible in $\Ho$ but not reversible in $\Hop$.

\subsection{Summary of results}\quad

The results in this paper are all stated for homeomorphisms of the real line. Most of the results remain valid when the domain is allowed to be any interval on the real line. 

Sections \ref{section:reversibility}, \ref{section:manyreversibility}, \ref{section:strongreversibility}, and \ref{section:involutions} are each structured similarly; the first subsection is about $\Ho$, the second subsection is about $\PLF(\R)$, and the third subsection is about $\PL(\R)$. Those results on $\Ho$ that are already known are stated and referenced, but not proved.

Section \ref{section:reversibility} is about reversible maps. In any group, a classification of reversible maps follows from a classification of conjugacy. It is known how to classify conjugacy in $\Ho$, $\PLF(\R)$, and $\PL(\R)$ (see \cite{BrSq01,Gl81}). We show in \S\ref{section:reversibility} that $\Hop$ and $\PL^+(\R)$ both contain non-trivial reversible maps whereas $\PLF^+(\R)$ contains none.

Section \ref{section:manyreversibility} is about the classes $R_n(G)$. Our findings are summarised in the next theorem.

\begin{theorem}\label{T: reversiblecomposites}
We have
\begin{enumerate}
\item $\Hop = R_2(\Hop)$ and $\Ho=R_2(\Ho)$;
\item for each integer $n$, $R_n(\PLF^+(\R))$ contains only the identity element, and  $I_{2n}(\PLF(\R))=R_n(\PLF(\R))$;
\item $\PL^+(\R)=R_4(\PL^+(\R))$ and $\PL(\R)=R_2(\PL(\R))$.
\end{enumerate}
\end{theorem}

Section \ref{section:strongreversibility} is about strong reversibility.
 The main results are summarised in the next theorem.  (Part (i) of
Theorem~\ref{T: strongReverse} follows from \cite[Theorem 2.2]{Yo94}.)

\begin{theorem}\label{T: strongReverse}\quad
\begin{enumerate}
\item In $\Ho$, an orientation preserving reversible map is strongly reversible if and only if it is reversible by an orientation reversing homeomorphism.
\item In $\PLF(\R)$, reversibility and strong reversibility are equivalent.
\item In $\PL(\R)$, an orientation preserving reversible map is strongly reversible if and only if it is reversible by an orientation reversing homeomorphism.
\end{enumerate}
\end{theorem}

In Section \ref{section:strongreversibility} we also prove the next theorem which is related to strong reversibility in $\Ho$.  Suppose that, given an element $f$ of a group $G$, one seeks involutions $\tau$ that are solutions of $\tau f\tau = f$, rather than $\tau f\tau = f^{-1}$. In a finite group this would be equivalent to determining whether the centraliser of $f$ is of even order. There is a neat characterisation of such maps $f$ in $\Ho$.

\begin{theorem}\label{T: squareRoot}
Given a map $f$ in $\Hop$, there is an involution $\tau$ in $\Hon$ such that $\tau f\tau =f$ if and only if $f$ has a square-root in $\Hon$.
\end{theorem}

Section \ref{section:involutions} is about the classes $I_n(G)$. A key result is summarised in the next theorem. (Part (i) of Theorem~\ref{T: composites} is the same as \cite[Theorem 6]{FiSc55}, \cite[Theorem 2]{Ja02a}, and \cite[Theorem 3.2]{OF04}.)

\begin{theorem}\label{T: composites}
We have
\begin{enumerate}
\item $\Hon \subsetneq I_3(\Ho)$ and $\Ho = I_4(\Ho)$,
\item $\PLF(\R)\supsetneq I_\infty(\PLF(\R))$,
\item $\PL^-(\R) \subsetneq I_3(\PL(\R))$ and $\PL(\R)=I_4(\PL(\R))$.
\end{enumerate}
\end{theorem}

\section{Reversibility}\label{section:reversibility}

We describe the reversible elements in our groups of interest.

\subsection{Reversibility in $\textup{H}$}\label{SS: reversible H}\quad

We first review necessary and sufficient conditions for functions in $\Ho$ to be conjugate. The \emph{degree}, $\textup{deg}_f$, of a homeomorphism $f$ is $1$ if $f$ preserves orientation, and $-1$ if $f$ reverses orientation. The \emph{signature} of $f$ is the function $\Gamma_f: \R\to\{-1,0,1\}$ given by
\[
\Gamma_f(x)=
\left\{
\begin{array}{ll}
1, & f(x)>x\,; \\
0, &  f(x)=x\,; \\
-1, & f(x)<x\,.
\end{array}
\right.
\]
Notice that $\Gamma_{f^{-1}}=-\Gamma_f$.  Conjugacy in both groups $\Ho$ and $\Hop$ can be understood in terms of the following two propositions (\cite[Lemmas 2.3 and 2.4]{OF04}).

\begin{proposition}\label{P: conjH+}
If $f$ and $g$ are two maps in $\Hop$, and $h$ is a member of $\Ho$ such that $hfh^{-1}=g$, then  $\Gamma_f = (\textup{deg}_h)\Gamma_g \circ h$. Conversely, if $h$ is a member of $\Ho$ such that $\Gamma_f = (\textup{deg}_h)\Gamma_g \circ h$, then $f$ and $g$ are conjugate in $\Ho$ by a homeomorphism of the same degree as $h$.
\end{proposition}

In particular, notice that the collection of fixed point free homeomorphisms consists of two conjugacy classes in $\Hop$---one containing the map $x\mapsto x+1$ and the other containing the map $x\mapsto x-1$---but these two conjugacy classes fuse in the larger group $\Ho$.  Similarly, the collection of fixed point free homeomorphisms in $\PL(\R)$ forms a single conjugacy class. In contrast, there are uncountably many conjugacy classes of fixed point free elements in $\PLF(\R)$. (See \cite{BrSq01} for proofs of these last two statements; the first is elementary.) 

\begin{proposition}\label{P: conjH-}
Two orientation reversing homeomorphisms $f$ and $g$ are conjugate in $\Ho$ if and only if there is a homeomorphism $h$ that maps the fixed point of $f$ to the fixed point of $g$ and is such that $hf^2h^{-1}=g^2$.
\end{proposition}

The only maps of finite order in $\Ho$, other than the identity, are orientation reversing involutions. By Proposition~\ref{P: conjH-}, all such maps are conjugate. 

We have two lemmas about the fixed point sets of maps $f$ and $h$ that satisfy the equation $hfh^{-1}=f^{-1}$. The first lemma is about all group actions. It is straightforward to prove.

\begin{lemma}\label{L: all}
Let $G$ be a group that acts on a set, and let $f$ and $h$ be elements of $G$ that satisfy $hfh^{-1}=f^{-1}$. Then $h$ permutes the fixed points of $f$.
\end{lemma}

Orientation preserving reversible homeomorphisms can be constructed using the next lemma.

\begin{lemma}\label{L: useful}
Suppose that $f,h\in\Hop$ and $hfh^{-1}=f^{-1}$. Then each fixed point of $h$ is also a fixed point of $f$.
\end{lemma}
\begin{proof}
Suppose that $h(p)=p$ for $p\in\mathbb{R}$. If $f(p)\neq p$ then, by swapping $f$ and $f^{-1}$ if necessary, we can assume that $f(p)>p$ and $f^{-1}(p)<p$. Therefore
\[
p>f^{-1}(p) = hfh^{-1}(p) = hf(p) > h(p) = p,
\]
which is a contradiction.  Hence $f(p)=p$.
\end{proof}

Given an orientation preserving homeomorphism $h$ we briefly describe a process for constructing another orientation preserving homeomorphism $f$ such that $hfh^{-1}=f^{-1}$. Each reversible homeomorphism in $\Hop$ arises by this process.  By Lemma~\ref{L: useful}, we must choose a map $f$ that fixes each   point of the fixed point set of $h$, $\text{fix}(h)$, in which case the  relationship $hfh^{-1}(x)=f^{-1}(x)$ certainly holds for points $x$ in $\text{fix}(h)$. Moreover, $f$ (like $h$) must fix, as a set,  each component $I$ in the complement of $\text{fix}(h)$. Now, $I$ is homeomorphic to $\mathbb{R}$, and, up to conjugacy, the only fixed point free orientation preserving homeomorphisms of $\mathbb{R}$ are $x\mapsto x+1$ and $x\mapsto x-1$. That is, there is an orientation preserving homeomorphism $k:I\rightarrow\mathbb{R}$ such that $khk^{-1}$ is either the map $x\mapsto x+1$ or the map $x\mapsto x-1$.  Suppose that we have found an orientation preserving homeomorphism $f_I$ of $I$ such that $hf_Ih^{-1}=f_I^{-1}$ on $I$.  Define $g=kfk^{-1}$. Then $g$ satisfies the functional equation $g(x+1)=g^{-1}(x)+1$ on $\mathbb{R}$. Conversely, given a solution $g$ of this functional equation, the function $f_I=k^{-1}gk$ satisfies the equation $hf_Ih^{-1}=f_I^{-1}$ on $I$. In summary, to solve the equation $hfh^{-1}= f^{-1}$ we must solve the functional equation $g(x+1)=g^{-1}(x)+1$ on each component of $\mathbb{R}\setminus\text{fix}(h)$.

We turn to reversibility in $\Ho$. According to Theorem~\ref{T: strongReverse} (i) (proven in \S\ref{SS: strongH}) a map in $\Hop$ that is reversible in $\Ho$ by an element of $\Hon$ is strongly reversible in $\textup{H}(\R)$. There are examples of orientation preserving homeomorphisms that are (i) not reversible;   (ii) reversible by an orientation preserving map but not reversible by an orientation reversing map;  (iii) reversible by an orientation reversing map but not reversible by an orientation preserving map;  or (iv) reversible by both orientation preserving and reversing maps. Here is a brief discussion of examples corresponding to these four possibilities. For (i), the map $f(x)=2x$ is not reversible, by Proposition~\ref{P: conjH+}, because if a homeomorphism $h$ reverses $f$ then it must fix $0$, in which case $\Gamma_{f^{-1}}\neq (\text{deg}_h) \Gamma_f \circ h$ . For (ii), first consider the map $f(x)=x+\sin(x)$. Let $t$ be the map $x\mapsto x+\pi$. Then $\Gamma_f\circ t = -\Gamma_f$. Since $-\Gamma_f=\Gamma_{f^{-1}}$ we see from Proposition~\ref{P: conjH+} that $f$ is reversible. Now consider any orientation preserving homeomorphism $k$ from $(0,+\infty)$ to $\mathbb{R}$. Define $g=kfk^{-1}$. Then $g$ is an orientation preserving homeomorphism of $(0,+\infty)$ which is reversed by another orientation preserving homeomorphism of $(0,+\infty)$ which we denote by $h$. Let us extend the definition of $g$ and $h$ to $\mathbb{R}$ by defining $g(x)=h(x)=x$ for each point $x\leq 0$. The equation $hgh^{-1}=g^{-1}$ is preserved. Suppose now that $sgs^{-1}=g^{-1}$ for an orientation \emph{reversing} homeomorphism $s$ of $\mathbb{R}$. There is a negative number $x_0$ such that for each number $x<x_0$ we have $s(x)>0$. For such points $x$ we obtain $s(x)=sg(x)=g^{-1}s(x)$. In other words, $g$ fixes each positive number greater than $s(x_0)$. This is false. Therefore $g$ is not reversible by an orientation reversing homeomorphism. For (iii), consider the map $x\mapsto \text{max}\left(2x,\tfrac{1}{2}x\right)$.  This map is reversed by $x\mapsto -x$ but, by Proposition~\ref{P: conjH+}, it is not reversible by an orientation preserving homeomorphism. For (iv), consider the map $f(x)= x+\cos x$. If $h$ denotes either of the maps $x\mapsto x+\pi$ or $x\mapsto -x$ then the equation $\Gamma_f=-(\text{deg}_h) \Gamma_f\circ h$ is satisfied.   From Proposition~\ref{P: conjH+} we deduce that $f$ is reversible by both orientation preserving and reversing homeomorphisms. (Of course, the identity map is also reversible by both orientation preserving and reversing homeomorphisms.)

Now consider $f\in\Hon$. If $hfh^{-1}=f^{-1}$ then also $(hf)f(hf)^{-1}=f^{-1}$, so $f$ is reversible by an orientation preserving map if and only if it is reversible by an orientation reversing map. By Proposition~\ref{P: conjH-},  $f$ is conjugate to $f^{-1}$ if and only if  $f^2$ and $f^{-2}$ are conjugate by a map that fixes the fixed point of $f$.
It is straightforward to construct examples of orientation reversing homeomorphisms that are either reversible or not reversible.

\subsection{Reversibility in $\PLF$}\quad

In \cite{BrSq01}, Brin and Squier characterise the conjugacy classes in $\PLF^+(\R)$ using geometrical invariants. In particular, for any function $f$ in $\PLF^+(\R)$, they describe  three quantities, which we collectively denote by $\Sigma_f$, with the property that maps $f$ and $g$ in $\PLF^+(\R)$ are conjugate if and only if $\Sigma_f=\Sigma_g$. One of these three quantities is $\Gamma_f$, which we have already encountered. The reader should consult \cite{BrSq01} for definitions of the other two quantities; they are not needed here.

It is sufficient to consider $\Gamma$ to understand reversibility in $\PLF^+(\R)$. For each map $f$ in $\PLF^+(\R)$, the function $\Gamma_f$ changes value only a finite number of times, so we can consider it to be a finite sequence with entries taken from $\{-1,0,1\}$. We exclude from this finite sequence those entries $0$ which correspond to isolated fixed points of $f$. (Thus an entry $0$ in the sequence corresponds to a closed interval of positive length on which $f$ coincides with the identity map.)  By Proposition~\ref{P: conjH+}, if two maps $f$ and $g$ in $\PLF^+(\R)$ are conjugate  then the two corresponding finite sequences must coincide. Since $\Gamma_f=-\Gamma_{f^{-1}}$, we see that the only reversible map in $\PLF^+(\R)$ is the identity map. 

There are, however, orientation preserving reversible maps in $\PLF(\R)$. For example, the map $x\mapsto \text{max}\left(2x,\frac12x\right)$ is reversed by $x\mapsto -x$. A complete classification of the reversible members of $\PLF(\R)$ follows immediately from the conjugacy classification in \cite{BrSq01}. 


Finally, we prove that the only orientation reversing reversible maps in $\PLF^-(\R)$ are involutions. In fact, we prove a stronger result that is valid for $\PL^-(\R)$. We use the notation $\text{fix}(g)$ to denote the fixed point set of a homeomorphism $g$, and we use the phrase \emph{bump domain} to describe a connected component of $\mathbb{R}\setminus\text{fix}(g)$.

\begin{proposition}\label{P: rev}
A member of $\textup{PL}^-(\R)$ is reversible in $\textup{PL}(\R)$ if and only if it is an involution.
\end{proposition}
\begin{proof}
Each involution in $\textup{PL}^-(\R)$ is reversed by the identity map. Conversely, suppose that $hfh^{-1}=f^{-1}$ for a map $f$ in $\PL^-(\R)$ and a map $h$ in $\textup{PL}(\R)$. By replacing $h$ with $hf$ if necessary, we can assume that $h$ preserves orientation. From the equation $hf(p) = f^{-1}h(p)$ we see that $h$ fixes the fixed point $p$ of $f$. Let $D$ be the first bump domain of $f^2$ to the right of $p$ on the real line, provided that such a domain exists. By Lemma~\ref{L: all} we know that $h$ maps one bump domain of $f$ to another. Since $h$ fixes $p$ and preserves orientation, $h$ must fix $D$ as a set. We arrive at a contradiction, since for a point $x$ in $D$, only one of $hf^2h^{-1}(x)$ and $f^{-2}(x)$ is greater than $x$. Similarly we can deduce that $f^2$ has no bump domains to the left of $p$ on the real line. Therefore $f^2$ is the identity map.
\end{proof}


\subsection{Reversibility in $\PL$}\quad

We make use of the following well known proposition. 

\begin{proposition}\label{P: freePL}
The fixed point free elements of $\PL^+(\R)$ fall into two conjugacy classes; one containing $x\mapsto x+1$ and the other containing $x\mapsto x-1$.
\end{proposition}

It is straightforward to prove Proposition~\ref{P: freePL} by constructing conjugating maps explicitly. Alternatively, in \cite{BrSq01} Brin and Squier state that the result can be proved using the techniques of that paper.

Suppose that $hfh^{-1}=f^{-1}$ for homeomorphisms $f$ and $h$ in $\PL^+(\R)$. Suppose that $h$ has a fixed point. Because $h$ is not the identity element, it must have a bump domain. Let us assume that this bump domain is of the form $(s,t)$, where $s$ is a real number and $t$ is either a real number or $+\infty$ (the case when $t\neq +\infty$, but possibly $s=-\infty$, is similar). By Lemma~\ref{L: useful}, $f$ fixes each point of $\text{fix}(h)$. Therefore $f$ fixes $(s,t)$. Suppose there is a point $p$ in $(s,t)$ such that $f(p)>p$. From \eqref{E: generalReverse} we see that $f(h^{n}(p))$ is greater than $p$ for even integers $n$ and $f(h^{n}(p))$ is less than $p$ for odd integers $n$. Since the points $h^n(p)$ accumulate at $s$, we have a contradiction, because members of $\PL^+(\R)$ are locally affine at all but a discrete set of points. Therefore $f$ is the identity map on $(s,t)$. Because we can repeat this procedure for any bump domain of $h$, and because, by Lemma~\ref{L: useful}, $f$ fixes each element of $\text{fix}(h)$, we deduce that $f$ is the identity map, in which case $h$ can be any element of $\PL^+(\R)$ whatsoever. 

The remaining possibility is that $h$ is free of fixed points. In this case, by Proposition~\ref{P: freePL}, either  $h$ or $h^{-1}$ is conjugate to $x\mapsto x+1$. Therefore $f$ is conjugate to a homeomorphism $g$ that satisfies $g(x+1)=g^{-1}(x)+1$ for each element $x$ of $\mathbb{R}$.

By Theorem~\ref{T: strongReverse} (iii), proven in \S\ref{SS: PLstrong}, an element of $\PL^+(\R)$ is reversed by an element of $\PL^-(\R)$ if and only if it is strongly reversible. Thus the orientation preserving reversible maps in $\PL(\R)$ are fully classified once \S\ref{SS: PLstrong} is complete. By Proposition~\ref{P: rev}, the only orientation reversing members of $\PL(\R)$ are the maps of order $2$. In the group $\PL(\R)$ it is straightforward to construct orientation preserving elements  that are (i) not reversible;   (ii) reversible by an orientation preserving map but not reversible by an orientation reversing map;  (iii) reversible by an orientation reversing map but not reversible by an orientation preserving map;  or (iv) reversible by both orientation preserving and reversing maps. We found such examples for the group $\Ho$ in \S\ref{SS: reversible H}, and similar examples can be constructed for $\PL(\R)$.

\section{Composites of reversible maps}\label{section:manyreversibility}

\subsection{Composites of reversible maps in $\textup{H}$}\quad

\begin{proof}[Proof of Theorem~\ref{T: reversiblecomposites} (i)]
First we prove that $\Hop=R_2(\Hop)$. Consider the two homeomorphisms $f$ and $g$ given by equations $f(x)=x+\tfrac{1}{4}$ and $g(x)=x+\tfrac{1}{2}\sin x$. Then $f$ is free of fixed points and both $g^{-1}$ and $fg$ have the same signature as the reversible map $x\mapsto x+\sin(x)$. By Proposition~\ref{P: conjH+} this means that $g^{-1}$ and $fg$ are both conjugate to $x\mapsto x+\sin(x)$, so they are both reversible. Hence $f$ lies in $R_2(\Hop)$. By conjugation we see that all fixed point free maps lie in $R_2(\Hop)$. Because open intervals are homeomorphic to \R, we see that each fixed point free homeomorphism $f_I$ of an open interval $I$ can be expressed as a composite of two reversible increasing homeomorphisms $g_I$ and $h_I$ of $I$.

Given $f$ in $\Hop$, let $\text{fix}(f)$ be the fixed point set of $f$, and define $f_I$ to be the restriction of $f$ to one of the open intervals $I$ in the complement of $\text{fix}(f)$. Define $g$ and $h$ in $\Hop$ as follows. For each $x$ in $\text{fix}(f)$, let $g(x)=h(x)=x$, and for each $x$ in $I$, let $g(x)=g_I(x)$ and $h(x)=h_I(x)$, where $g_I$ and $h_I$ are the functions described in the last sentence of the preceding paragraph. Then $f=gh$.

That $\Ho=R_2(\Ho)$ follows from Theorem~\ref{T: composites} (i) (proved in \S\ref{SS: compositesHR}), because $R_2(\Ho)$ contains $I_4(\Ho)$.
\end{proof}

\subsection{Composites of reversible maps in $\PLF$}\quad

There are no reversible homeomorphisms in $\PLF^+(\R)$ other than the identity. According to Theorem~\ref{T: strongReverse} (ii) each reversible homeomorphism in $\PLF(\R)$ is strongly reversible. Therefore the sets  $R_n(\PLF(\R))$ coincide with $I_{2n}(\PLF(\R))$ for each integer $n$. We examine such sets in \S\ref{SS:srPLF} and \S\ref{SS:involutionsPLF}. (The content of Theorem~\ref{T: reversiblecomposites} (ii) has now been addressed.)

\subsection{Composites of reversible maps in $\PL$}\quad

\begin{proposition}\label{P: PLcomposites}
Each fixed point free homeomorphism in $\PL^+(\R)$ can be expressed as a composite of two reversible homeomorphisms in $\PL^+(\R)$.
\end{proposition}
\begin{proof}
Let $f$ the unique member of $\PL^+(\R)$ such that
\[
f(x)=
\begin{cases}
3x, & 0\leq x\leq 2\,;\\
\tfrac{1}{3}(x+16), & 2\leq x\leq 8\,;
\end{cases}
\]
and for all real numbers $x$ the equation $f(x+8)=f^{-1}(x)+8$ is satisfied. The map $f$ has been constructed so that it is reversed by the translation $x\mapsto x+8$. Define a second member of $\PL^+(\R)$ by the equation $g(x)=f(x)-2$. This map $g$ is conjugate to $f$ in $\PL^+(\R)$. One may see this by appealing to \cite{BrSq01}, or alternatively one can construct the conjugating map explicitly: it is the unique member of $\PL^+(\R)$ which, for each integer $n$,  is locally affine on  $(16n+1,16n+5)$ and maps this interval to  $(16n,16n+8)$, and is locally affine on $(16n+5, 16(n+1)+1)$ and maps this interval to $(16n+8, 16(n+1))$. Since $g$ is conjugate to $f$, it is also  reversible.

We have $g^{-1}f(x)>x$ for each $x$.  Thus, by Proposition~\ref{P: freePL}, a fixed point free homeomorphism in $\PL^+(\R)$ is conjugate to  either $g^{-1}f$ or $f^{-1}g$, and as such it can be expressed as a composite of two reversible homeomorphisms in $\PL^+(\R)$.
\end{proof}

\begin{proof}[Proof of Theorem~\ref{T: reversiblecomposites} (iii)]
First we show that $\PL^+(\R) = R_4(\PL^+(\R))$. Choose $f$ in $\PL^+(\R)$. Choose a fixed point free map $g$ in $\PL^+(\R)$ such that $f(x)>g(x)$ for each real number $x$. Then $f=gh$, where $h=g^{-1}f$, and both $g$ and $h$ are fixed point free members of $\PL^+(\R)$. The result follows from Proposition~\ref{P: PLcomposites}.

Second we show that $\PL(\R) = R_2(\PL(\R))$. This follows from Theorem~\ref{T: composites} (iii), proven in \S\ref{SS: involutionsPL}, and the fact that $R_2(\PL(\R))$ contains $I_4(\PL(\R))$.
\end{proof}

We do not have a characterisation of the sets $R_2(\PL^+(\R))$ and $R_3(\PL^+(\R))$.

\section{Strongly reversible maps}\label{section:strongreversibility}

There are no involutions in $\textup{H}^+(\R)$ other than the identity thus, in order to make worthwhile statements in this section,  we work in the full homeomorphism groups $\Ho$, $\PLF(\R)$, and $\PL(\R)$. Furthermore, since the composite of two non-trivial involutions in $\textup{H}(\R)$ is a map in $\textup{H}^+(\R)$, the only strongly reversible orientation reversing homeomorphisms in $\Ho$, $\PLF(\R)$, and $\PL(\R)$ are involutions. Thus we are left to examine the orientation preserving strongly reversible members of $\Ho$, $\PLF(\R)$, and $\PL(\R)$.

\subsection{Strong reversibility in $\textup{H}$}\label{SS: strongH}\quad

The orientation preserving strongly reversible maps in \Ho can be classified according to the following proposition, which is an extension of  \cite[Theorem 1]{Ja02b} and \cite[Theorem 2.2]{Yo94}.

\begin{proposition}\label{P: strongrev h+}
For each map $f$ in $\textup{H}^+(\R)$, the following are equivalent:
\begin{enumerate}
\item $f$ is strongly reversible;
\item there is an orientation reversing homeomorphism $h$ such that $\Gamma_f = \Gamma_f\circ h$;
\item $f$ is conjugate to a homeomorphism that is strongly reversed by the involution $\eta(x)=-x$;
\item $f$ is conjugate to a homeomorphism whose graph is symmetric across the line $y=-x$.
\end{enumerate}
\end{proposition}
\begin{proof}
The equivalence of (i) and (ii) is proven in \cite[Theorem 2.2]{Yo94}. Statements (i) and (iii) are equivalent because, by Proposition~\ref{P: conjH-}, all non-trivial involutions are conjugate. It remains only to show that (iii) and (iv) are equivalent. Let $G(k)$ denote the subset of $\mathbb{R}^2$ that is the graph of a homeomorphism $k$. Notice that $G(\eta f\eta)$ is $G(f)$ rotated by an angle $\pi$ about $0$, and $G(f^{-1})$ is the reflection of $G(f)$ in the line $y=x$. That (iii) and (iv) are equivalent follows because $\eta f\eta =f^{-1}$ if and only if $G(\eta f\eta)=G(f^{-1})$.
\end{proof}

Conditions (iii) and (iv) add little substance to Jarczyk and Young's theorem (\cite[Theorem 1]{Ja02b} and \cite[Theorem 2.2]{Yo94}), but we include them as they are simpler to grasp than (ii), and they highlight a geometric connection between involutions (homeomorphisms whose graphs are symmetric across the line $y=x$) and strongly reversible maps (homeomorphisms that are \emph{conjugate} to maps whose graphs are symmetric across the line $y=-x$).

\begin{proof}[Proof of Theorem~\ref{T: strongReverse} (i)]
Suppose that $f\in\Hop$ is reversed by $h\in\Hon$. Observe that $\Gamma_f=-\Gamma_{f^{-1}}$. Using Proposition~\ref{P: conjH+} we see that
\[
\Gamma_f=-\Gamma_{f^{-1}} = \Gamma_f\circ h^{-1}.
\]
Therefore $f$ is strongly reversible, by Proposition~\ref{P: strongrev h+}.
\end{proof}

We finish this section with a proof of Theorem~\ref{T: squareRoot}.

\begin{proof}[Proof of Theorem~\ref{T: squareRoot}]
Suppose that $\tau f\tau=f$ for some involution $\tau\in\Hon$. Because $f$ is orientation preserving, it can be embedded in a flow (see \cite[Theorem 1]{FiSc55}), hence there is a map $h$ in $\Hop$ such that $h^2=f$. Let $p$ be the fixed point of $\tau$. Then $p$ is also a fixed point of $f$ and $h$. Define a map $g$ in $\Hon$ by the equation
\[
g(x)=
\begin{cases}
h\tau(x), & x\geq p\,;\\
\tau h(x), & x<p\,;
\end{cases}
\]
so that $g^2=f$. Conversely, if $g\in\Hon$ and $g^2=f$ then define $\tau\in\Hon$ by the equation
\[
\tau(x)=
\begin{cases}
g(x), & x\geq p\,;\\
g^{-1}(x), & x<p\,;
\end{cases}
\]
where $p$ is the fixed point of $g$. Then $\tau$ is an involution and $\tau f\tau =f$.
\end{proof}

\subsection{Strong reversibility in $\PLF$}\label{SS:srPLF}\quad

In this section we prove Theorem~\ref{T: strongReverse} (ii). For $f\in\PLF^+(\R)$ recall that the \emph{bump domains} of $f$ are the connected components of $\mathbb{R}\setminus\text{fix}({f})$. Furthermore we say that $f$ is a \emph{one bump function} if $\mathbb{R}\setminus\text{fix}({f})$ is connected and non-empty.

\begin{lemma}\label{L: singleBump}
Suppose that $f$ is a one bump function with bump domain $I$, and $h$ is a member of $\PL^-(\R)$  such that $hfh^{-1}=f^{-1}$. Then $h$ fixes $I$, and acts as an involution on $I$.
\end{lemma}
\begin{proof}
The bump domain $I$ is the complement of $\text{fix}(f)$. Hence, by Lemma~\ref{L: all}, it is fixed, as a set, by $h$. Since $h$ is orientation reversing, it has a fixed point in $I$.

Now, the map $h^2$ fixes $I$ and it commutes with $f$, by \eqref{E: generalReverse}. According to \cite[Theorem 4.18]{BrSq01}, the centraliser of $f$ in the group of piecewise linear homeomorphisms of $I$  is a cyclic group. Therefore there are integers $m$ and $n$ such that $f^m=h^{2n}$ on $I$. But $h^2$ has a fixed point in $I$, whereas $f^m$ does not, for each $m\neq 0$. Hence $h^2$ is the identity map on $I$. 
\end{proof}

\begin{proof}[Proof of Theorem~\ref{T: strongReverse} (ii)]
We have to show that if a map $f$ in $\PLF(\R)$ is reversible  in $\PLF(\R)$ then it is strongly reversible. If $f$ is orientation reversing then, by Proposition~\ref{P: rev}, it is an involution, which means that it is strongly reversible. Suppose then that $f\in\PLF^+(\R)$. Only the identity map is reversible in $\PLF^+(\R)$, so we assume that $hfh^{-1}=f^{-1}$, for $h\in\PLF^-(\R)$. Let $h$ have fixed point $p$. Define an involution $\tau$ in $\PLF^-(\R)$ by the equation
\begin{equation}\label{E: tauPL}
\tau(x)=
\begin{cases}
h^{-1}(x), & x\geq p\,;\\
h(x), & x<p\,.
\end{cases}
\end{equation}
We have only to show that $\tau f\tau =f^{-1}$. If $p$ is a fixed point of $f$ then certainly $\tau f\tau=f^{-1}$. If $p$ is not a fixed point of $f$  then $p$ lies inside a bump domain $(s,t)$ of $f$. Since, by Lemma~\ref{L: all}, $h$ fixes $\text{fix}(f)$ setwise, we see that $h$ fixes $(s,t)$ setwise. From Lemma~\ref{L: singleBump} we deduce that $h(x)=h^{-1}(x)$ for points $x$ in $I$, so the equation $\tau f\tau = f^{-1}$ is still satisfied.
\end{proof}

\subsection{Strong reversibility in $\PL$}\label{SS: PLstrong}\quad

\begin{proof}[Proof of Theorem~\ref{T: strongReverse} (iii)]
We need only show that if a map $f$ in $\PL^+(\R)$ is reversible in $\PL(\R)$ by an element $h$ of $\PL^-(\R)$, then $f$ is strongly reversible. Define $\tau$ by \eqref{E: tauPL}. Then identical reasoning to that given in the proof of Theorem~\ref{T: strongReverse} (ii) can be used to show that $\tau f\tau = f^{-1}$.
\end{proof}

\section{Composites of Involutions}\label{section:involutions}

\subsection{Composites of involutions in $\textup{H}$}\label{SS: compositesHR}\quad

Theorem~\ref{T: composites} (i) is equivalent to \cite[Theorem 6]{FiSc55} and \cite[Theorem 2]{Ja02b}. The proof, which we omit, is very similar to our proof of Theorem~\ref{T: composites} (iii), in \S\ref{SS: involutionsPL}. 

\subsection{Composites of involutions in $\PLF$}\label{SS:involutionsPLF}\quad

The results in this section differ substantially from those for  $\textup{H}(\R)$. We define two new subgroups of $\textup{PLF}^+(\R)$. Let $\textup{PLc}^+(\R)$ denote the group consisting of those members of $\PLF^+(\R)$ which coincide with the identity map outside a compact subset of $\mathbb{R}$. Next, we denote the left-most and right-most gradients of a map $f$  in $\PLF(\R)$ by $\gamma_L(f)$ and $\gamma_R(f)$. Define
\[
\PLE(\R)=\left\{f\in \PLF(\R)\,\big|\, \gamma_R(f)=\frac{1}{\gamma_L(f)}\right\}.
\]

\begin{proposition}\label{P: infinity}
We have $I_\infty(\PLF(\R))=\PLE(\R).$
\end{proposition}
\begin{proof}
First we prove that $I_\infty(\PLF(\R))\subseteq \PLE(\R)$. Let $\tau_1,\dots,\tau_m$ be involutions in $\PLF(\R)$. For each $i$, $\gamma_R(\tau_i)=\gamma_L(\tau_i)^{-1}$, therefore
$\gamma_L(\tau_1\cdots \tau_m)=\gamma_R(\tau_1\cdots \tau_m)^{-1}$ by the chain-rule.

Second we prove that $\PLE(\R)\subseteq I_\infty(\PLF(\R))$. Choose $f$ in $\PLE(\R)$. By composing with $\eta(x)=-x$ if necessary, we assume that $f$ preserves orientation. For each $p>0$, define an involution $\tau_p$ by the equation
\[
\tau_p(x) =
\begin{cases}
-p x, & x\geq 0\,;\\
-\tfrac{1}{p}x, & x<0\,.
\end{cases}
\]
For each real number $t$, define an involution $\sigma_t$ by the equation $\sigma_t(x)=-x+t$. For sufficiently large real numbers $x$ there are constants $\lambda>0$ and $u\in\mathbb{R}$ such that $f(x)=\lambda x  + u$. Likewise, for sufficiently small real numbers $x$ there is a constant $v\in\mathbb{R}$ such that $f(x)=\frac{1}{\lambda}x + v$. Define $p$ to be the unique positive number such that $p-\frac1p = u-v$. Define $q=\frac{\lambda}{p}$, and define $s=u-p$. Finally, define $g=\sigma_s\tau_p\sigma_1\tau_q$. One can check that for sufficiently large and small values of $x$ in $\R$, we have $f(x)=g(x)$. Therefore $fg^{-1}\in\PLc^+(\R)$. 

It remains to show that $\PLc^+(\R)$ is a subset of $I_\infty(\PLF(\R))$.  According to \cite[Theorem 3.1]{Ep70}, the group $\PLc^+(\R)$ is simple.  But $\PLc^+(\R)\cap I_\infty(\PLF(\R))$ is a normal subgroup of $\PLc^+(\R)$, and one can check that this subgroup contains the element $k=\sigma_{-5}\tau_2\sigma_2\tau_{1/4}\sigma_2\tau_2$. (To perform this check, first verify that $k(-1)=-1/4$---so that $k$ is not the identity map---then show that $k(x)=x$ for sufficiently large and small values of $x$.) Therefore 
$\PLc^+(\R)\cap I_\infty(\PLF(\R))=\PLc^+(\R)$, as required.
\end{proof}

Because $\PLE(\R)$ is strictly contained in $\PLF(\R)$,  Theorem~\ref{T: composites} (ii) follows from Proposition~\ref{P: infinity}.

\subsection{Composites of involutions in $\PL$}\label{SS: involutionsPL}\quad

\begin{proof}[Proof of Theorem~\ref{T: composites} (iii)]
Given an element $f$ of $\PL^-(\R)$, choose an involution $\sigma$ in $\PL^-(\R)$ such that, for each real number $x$, $\sigma(x) > f(x)$. Therefore $\sigma f(x)>x$. By Proposition~\ref{P: freePL}, $\sigma f$ is conjugate to the translation $x\mapsto x+1$.  This translation is reversed by the involution $x\mapsto -x$. Therefore $\sigma f\in I_2(\PL(\R))$, which means that $f\in I_3(\PL(\R))$. Both parts of the result follow immediately.
\end{proof}

Since all fixed-point free elements of $\PL(\R)$ are conjugate, the method of the proof of Theorem~\ref{T: composites} (i) given in \S\ref{SS: compositesHR} can also be used to prove Theorem~\ref{T: composites} (iii). 

\section{Open questions}

Two questions from our study remain open. In Theorem~\ref{T: reversiblecomposites} (iii) we showed that
$\PL^+(\R)=R_4(\PL^+(\R))$. 

\begin{opening}
What is the smallest integer $n$ for which $\PL^+(\R)=R_n(\PL^+(\R))$?
\end{opening}

In Proposition~\ref{P: infinity} we showed that $\textup{PLE}(\R) = I_\infty(\textup{PLF}(\R))$.

\begin{opening}
Is there an integer $n$ for which $\textup{PLE}(\R) = I_n(\textup{PLF}(\R))$?
\end{opening}

There are other related groups of homeomorphisms of the line and circle, such as piecewise linear homeomorphisms of the circle and Thompson's groups, for which questions of reversibility are open. A document on reversibility in the group of homeomorphisms of the circle is in preparation (\cite{GSO}).


\bibliographystyle{plain}
\bibliography{reversible}

\end{document}